\documentclass[14pt]{article}
\usepackage{graphicx}
\usepackage{url}
\usepackage[russian,english]{babel}
\usepackage{amsmath,latexsym,amsthm}
\usepackage{inputenc}
\usepackage{graphics}
\usepackage{graphicx}
\usepackage{amssymb}
\usepackage{color}
\usepackage[colorlinks]{hyperref}
\newtheorem{theorem}{Theorem}[section]
 \newtheorem{definition}[theorem]{Definition}

\newtheorem{lemma}[theorem]{Lemma}
\newtheorem{corollary}[theorem]{Corollary}
\newtheorem{remark}[theorem]{Remark}
\newtheorem{example}[theorem]{Example}
\begin{document}
\newcommand{\tit}[1]{
\begin{center} { {\bf #1 }}  \end{center}}  
\newcommand{\aut}[1]{
\begin{center} {\bf #1 }  \end{center}}           
\tit{ON LIMIT THEOREMS FOR FUNCTIONAL AUTOREGRESSIVE PROCESSES WITH RANDOM COEFFICIENTS}
\aut{Sadillo Sharipov}
\begin{center}
\emph{V.I.Romanovskiy Institute of Mathematics, Tashkent, Uzbekistan}
\end{center}
\begin{center}
E-mail: sadi.sharipov@yahoo.com
\end{center}

\newcommand{\Abstract}
\noindent \textbf{\Abstract}{ \small Abstract.}
In this paper, we consider a Banach space valued random coefficient autoregressive process. Our  studies on this process involve existence, weak law of large numbers, strong law of large numbers, some exponential inequalities, central limit theorem. Our approach is based on a suitable martingale coboundary decomposition in Banach space.\\

\newcommand{\keywords}
\\

\section{Introduction}
In recent years Functional Data Analysis (FDA) has established itself as an important field of modern statistics because of its applicability to problems which are difficult to cast into a framework of scalar or vector observations. Functional data come in many forms, however it always consist of functions, often smooth curves. Quite regularly, functional data are collected sequentially over time, and different curves of functional observations time record into disjoint natural time intervals. The monograph of Ramsey and Silverman \cite{book 5} has become a classic reference to the ideas and a large review of functional techniques of FDA. Other key references include  monographs by Bosq \cite{book 1}, Bosq and Blanke \cite{book 2}, Ferraty and Vieu \cite{paper 19}, Horv\'{a}th and Kokoszka \cite{paper 25} and surveys Cuevas \cite{paper 14}, Goia \cite{paper 20}, Gonz\'{a}lez-Manteiga \cite{paper 21} and reference therein.

As quoted below, since a variety of functional data is collected sequentially over time then we expect that the functional data in a given time period are effected by past observations. Then as the main tool of analyze data is used models of Functional Time Series (FTS). The literature on FTS has mainly centered around stationary linear models Bosq \cite{book 1}, Bosq and Blanke \cite{book 2}.

The simplest model for a FTS is Functional Autoregressive Processes with deterministic operator. This class extends to the functional setting the autoregressive model AR(1) and it is a very flexible modeling and predictive tool for continuous time random process.

The class of Functional Autoregressive Processes with deterministic operator is a very flexible modeling and predictive tool for continuous time random process. The general theory of this model was presented in the pioneering work of Bosq \cite{book 1} where developed estimation of its second order structure and derived an asymptotic theory. This model was successfully applied in road traffic \cite{paper 8}, climatology \cite{paper 7}, and to predict electricity consumption \cite{paper 13}.

This paper deals with Functional Autoregressive Processes with random coefficients. The random coefficient time series have been used in context of random perturbations of dynamical systems and they found a numerical of applications, such as economy, finance and biology etc. (see Nicholls and Quinn \cite{book 4}, Tj{\o}steim \cite{paper 35}, Tong \cite{book 6} and references therein). For Functional Autoregressive Processes with random coefficient Guillas \cite{paper 23} considered a model for which the random coefficients take two operator values with Bernoulli distribution and obtained some limit theorems for that model. Later, Mourid \cite{paper 30} studied a more general model for which the sequence of random coefficients is reduced to one random operator. Cugliary \cite{paper 15} extended the model of Guillas \cite{paper 23} to the case with many regimes. Recently, Boukhair and Mourid \cite{paper 10} considered a Hilbert space valued random coefficient autoregressive process and obtained some limit theorems for this process. However, their results require some additional restrictive assumptions such as boundedness of process and condition \eqref{eq:condition of Mourid} (see below).

The main purpose of this paper is to derive some limit theorems for Banach-valued autoregressive processes with random coefficients. Our results improve results of Boukhiar and Mourid \cite{paper 10} in two aspects. First, we impose weaker condition on random coefficients than Boukhair and Mourid \cite{paper 10} and at the second, we do not require on process to be bounded.

Our approach is based on a coboundary decomposition of Banach valued random coefficient autoregressive process.
Basicly, the main idea is to use martingale approximation for Banach-valued random variables. There have been many research works on martingale approximation for stationary processes. Early contribution is due to Gordin \cite{paper 22}, who proved CLT for stationary process using martingale approximation. In the literature, this method is also known as Gordin's method. We refer for further contributions to Hannan \cite{paper 24}, Maxwell-Woodroofe \cite{paper 27}, Peligrad and Utev \cite{paper 32}, Wu and Woodroofe \cite{paper 39}, Voln\'{y} \cite{paper 36} and among others where sharp results obtained concerning CLT and other limit theorems. For Banach space setting, this method was successfully applied to derive limit theorems; see e.g. Merlev\'{e}de \cite{paper 29}, Dedecker \cite{paper 17}, Dedecker and Merlev\'{e}de \cite{paper 18} and reference therein. A key reference is a recent monograph by Merlev\'{e}de et al. \cite{paper 28} where one can find recent developments of martingale approximation.

The paper is organized as follows. In Section 2 we provide basic facts and definitions. Section 3 contains main results and in Section 4 we provide proofs of the main results.

\section{Notations and definitions}

We begin by introducing basic facts. Let $(\Omega,\mathcal{A},\textbf{\textrm{P}})$  be a complete probability space and  $\textbf{B}$ be a separable Banach space with norm $\|\cdot\|_{\textbf{B}}$  and let $\textbf{B}^{*}$ denote its dual space. For any real $p\geq1$  denote
$ \textbf{\textrm{L}}_\textbf{B}^{p}$ the space of $\textbf{B}-$valued random variables such that $\|X\|_{{\textbf{\textrm{L}}}_{\textbf{B}}^{p}}^{p}=\textbf{E}\|X\\|_{\textbf{B}}^{p}$  is finite. We denote by $\mathcal{L}(\textbf{B})$  the space of bounded linear operators over $\textbf{B}$ equipped with usual uniform norm.

Let $(\rho_{n},n\in\mathbb{Z})$  be a sequence of measurable random operators defined on $(\Omega,\mathcal{A},\textbf{\textrm{P}})$ with values in $\mathcal{L}(\textbf{B})$  endowed with its Borel $\sigma-$field.

We first define a Banach-valued white noise.

\begin{definition}
A sequence $(\varepsilon_{n},n\in\mathbb{Z})$ of $ \textbf{B}-$valued random variables is called a strong  white noise (innovation process) if it is an i.i.d. and such that $\textbf{E}\varepsilon_{0}=0$ and $\textbf{E}\|\varepsilon_{0}\|_{\textbf{B}}^{2}<\infty$.
\end{definition}
Now recall the definition of Banach-valued random coefficients autoregressive processes of order 1.
\begin {definition}
A sequence $(X_{n},n\in\mathbb{Z})$ of $\textbf{B}-$valued random variables satisfying the following recursion equation
\begin{equation}
X_{n}-\mu=\rho_{n}(X_{n-1}-\mu)+\varepsilon_{n}, \ \ n\in\mathbb{Z}\label{eq:autoregressive equation}
\end{equation}		
where $\mu\in \textbf{B}$, is called a Banach-valued random coefficients autoregressive  process of order 1 (we abbreviate by BRCA(1)).
\end{definition}
The sequence \eqref{eq:autoregressive equation} defines a nonlinear functional time series model, thus extends the class of functional autoregressive models paramount in linear functional time series analysis. The process defined by \eqref{eq:autoregressive equation} is used to handle possible nonlinear features of data and includes nonlinear models, threshold models and some double stochastic time series.

We recall some known facts.
\begin{definition}
For $\textbf{B}-$valued random variable $X$ with $\textbf{E}X=0 $ and $\textbf{E}\|X\|_{\textbf{B}}^{2}<\infty$, covariance operator $C_{X}$ of $X$ is the bounded linear operator from $\textbf{B}^{*}$ to $\textbf{B}$, defined by $C_{X}(x^{*})=\textbf{E}(x^{*}(X)X)$,
$x^{*} \in \textbf{B}^{*}$.
\end{definition}
Let $ \textbf{B}_{1}$  and $\textbf{B}_{2}$ be two separable Banach spaces, let $X$ and $Y$ be random variables with zero means and belong to $\textbf{\textrm{L}}_{\textbf{B}_{1}}^{2}$  and $\textbf{\textrm{L}}_{\textbf{B}_{2}}^{2}$, respectively. Then the cross-covariance operators of $X$ and $Y$ are defined by  $C_{X,Y}(x^{*})=\textbf{E}(x^{*}(X)Y)$, $ x^{*} \in \textbf{B}_{1}^{*}$.
We use $\left(\mathcal{L}^\star,\|\cdot \|^\star \right)$ to denote the space of bounded linear operators from $\textbf{B}^{*}$  to $\textbf{B}$ equipped with uniform norm
$$ \|\rho\|^\star=\sup_{\|x^{\ast}\|\leq 1}{\|\rho^\star(x^{\ast})\|},\ \ \rho^\star\in \mathcal{L}^\star. $$
An operator $\rho\in \mathcal{L}^\star $ is called nuclear if it admits the representation
$$ \rho\left(x^\ast\right)=\sum_{j=1}^{\infty}{x_{j}^{\ast\ast}\left(x^{\ast }\right)y_{k}},\ \ x\in B^\ast, $$
where $\left(x_{k}^{\ast\ast} \right)\subset B^{\ast\ast}-$ the dual space of $B^{\ast}$, $\left(y_{k}\right)\subset B$ with $\sum_{j=1}^{\infty}{\|x_{j}^{\ast\ast}\|\|y_{k}\|}<\infty.$ The infimum of sums $\sum_{j=1}^{\infty}{\|x_{j}^{\ast\ast}\|\|y_{k}\|}$ for all such representation of $\rho $ is a norm. We use $(\mathcal{N},\|\cdot\|_{\mathcal{N}})$ to denote the set of nuclear operators with norm $\|\cdot\|_{\mathcal{N}}$.

Covariance operator is symmetric, compact, positive, and nuclear.
\begin{definition}
Let $(\mathcal{M}_{i},i\geq 1)$  be an increasing sequence of sub-$ \sigma-$algebras of $\mathcal{A}$. We say that $(\xi_{i}, i\geq 1)$  is a sequence of $ \textbf{B}-$valued martingale differences with respect to filtration $(\mathcal{M}_{i},i\geq 0)$ if

1. for all $i\geq1$, $\xi_{i}$ is $ \mathcal{M}_{i}$ measurable and belongs to $\textbf{\textrm{L}}_\textbf{B}^{1}$;

2. for all  $i>1$, $\textbf{E}(\xi_{i}|\mathcal{M}_{i-1})=0$  almost surely.
\end{definition}

We recall the notation of $r-$smooth Banach spaces, $1<r\leq2$.

\begin{definition}\label{dfn:r_smooth}
We say that a separable Banach space $(\textbf{B}, \|\cdot\|_{\textbf{B}})$ is $r-$smooth $(1<r\leq2)$ if there exists an equivalent norm $\|\cdot\|^{'}$  such that
\begin{center}
$\mathop {\sup }\limits_{t>0}\{ \frac{1}{t^{r}}{\sup}\{{||x+ty||^{'}+||x-ty||^{'}-2: ||x||^{'}=||y||^{'}=1}\}\}<\infty.$
\end{center}
\end{definition}
This notion was introduced by Pisier \cite{paper 33}. Note that if $ \textbf{B}$ is $r-$smooth then it is $r^{\prime}-$smooth for any $r^{\prime}<r$. According to \cite{paper 5}, if  $\textbf{B}$  is $r-$smooth and separable then there exist a constant  $D\geq1$ such that, for any sequence of $ \textbf{B}-$valued martingale difference $(\xi_{i}, i\geq 1)$
\begin{center}
$ \textbf{E}||\xi_{1}+...+\xi_{n}||_{\textbf{B}}^{r}\leq D \sum \limits_{i = 1}^n {\textbf{E}||\xi_{i}||_{\textbf{B}}^{r} }$.
\end{center}
Therefore, we may claim that these spaces play the same role with respect to martingales as spaces to type $p$  do with respect to the sums of independent variables.


\section{Main results}

The basic idea of martingale approximation of sums of stationary process is following. Let $T:\Omega\rightarrow\Omega $ be a bijective, bimeasurable and measure-preserving transformation, $f$ is a stationary process that belongs to $\textbf{L}^{p}, p\geq1$ and consider the random sequence $(f\circ T^{n})_{n\in\mathbb{Z}}$. If $(f\circ T^{n})_{n\in\mathbb{Z}}$ is a $p-$integrable stationary process then $f$ admits a coboundary decomposition:
\begin{equation}\label{eq:martingale_coboundary}
f=h+g-g\circ T
\end{equation}
where $g$ is a $p-$integrable function and $(h\circ T^{n})_{n\in\mathbb{Z}}$ is a $p-$integrable stationary martingale difference sequence. The term $g-g\circ T$ is called a coboundary.
The decomposition \eqref{eq:martingale_coboundary} allows to prove limit theorems for $f$ from corresponding results for stationary martingale difference sequences. There exist necessary and sufficient condition for $f$ to be fulfilled a coboundary decomposition \eqref{eq:martingale_coboundary}; see e.g. \cite{paper 26}. For instance, when $p=2$ this criterion is the same as that given by Gordin \cite{paper 22} for CLT.

Now from suitable coboundary decomposition of BRCA(1) process, we derive martingale approximation for the sample mean of \eqref{eq:autoregressive equation} that will be the basis for all further analysis.
For each $i\geq1$, let $\mathfrak{B}_{i}$  be an increasing filtration of $\sigma-$algebras generated by $((\rho_{j},\varepsilon_{j}),j\leq i)$.

For $(X_{n},n\in\mathbb{Z})$ process defined by equation \eqref{eq:autoregressive equation} consider the decomposition
\begin{equation}
X_{i}=N_{i-1}-N_{i}+M_{i},\ \ i\geq1,\label{eq:autoregressive decomposition}
\end{equation}
where
\begin{center}
$N_{i}=(I-\overline{\rho})^{-1}\overline{\rho}X_{i}, \ \   M_{i}=(I-\overline{\rho})^{-1}X_{i}-(I-\overline{\rho})^{-1}\overline{\rho}X_{i-1},$
\end{center}
with $\overline{\rho}=\textbf{E}(\rho_{0})$, $||\overline{\rho}||_{\mathcal{L}(B)}\leq \textbf{E}||\overline{\rho}||_{\mathcal{L}}<1$  and thus $(I-\overline{\rho})^{-1}$ is well defined over $ \textbf{B}$. Note that since $ \varepsilon_{i}$ is an independent of $ \mathfrak{B}_{i-1}$  and random variable $X_{i-1}$ is a measurable with respect to $ \mathfrak{B}_{i-1}$, we have  $\textbf{E}(X_{i}|\mathfrak{B}_{i-1})=\overline{\rho}X_{i-1}$. This implies that the variables $(M_{i}, i\geq 1)$  defines a sequence of $ \textbf{B}-$valued stationary martingale differences with respect to the filtration $\mathfrak{B}_{i-1}$.
Adding the identities in \eqref{eq:autoregressive decomposition}, we get
\begin{equation}
\sum \limits_{i=1}^n{X_{i}}=N_{0}-N_{n}+\sum \limits_{i=1}^n{M_{i}}\label{eq:sum of autoregressive process}.
\end{equation}				
Hence, the sum of $X_{i}$ admits a suitable martingale approximation and we may apply recent results on martingale approximation to \eqref{eq:sum of autoregressive process}.

In the first part of this paper, we consider estimation of the mean $\mu$ from observations $X_{1},...,X_{n}$, where $(X_{n})$ is BRCA(1) process. The natural estimator is the sample mean
\begin{center}
$ {\overline{X}_{n}}:=\frac{1}{n}\sum\limits_{i=1}^{n}{X_{i}}.$
\end{center}
Observe that $\textbf{E}\overline{X}_{n}=\mu $, that is, that $\textbf{E}\overline{X}_{n}$ is an unbiased estimator of $\mu $. Throughout the paper, we set $S_{n}:=X_{1}+...+X_{n}$, $n\geq1$.

We now give our main results. We first study sufficient conditions for the existence and uniqueness of strictly stationary BRCA(1) process. In order to provide the existence and uniqueness of strictly stationary solution of BRCA(1) we need the following conditions:
\begin{enumerate}
\item[(C1($p$)):]\label{cond:independence_of_rho} The random variables $(\rho_{n},n\in\mathbb{Z})$  are i.i.d. and belongs to $\textbf{\textrm{L}}_\textbf{B}^{p}$, $p\geq1$.
\item[(C2):] The two sequences $(\rho_{n},n\in\mathbb{Z})$  and $(\varepsilon_{n},n\in\mathbb{Z})$  are independent.
\item[(C3($p$)):] \label{cond:cond_on_p_moment_of_rho} $\textbf{E}||\rho_{0}||_\mathcal{L}^{p}<1$, $p\geq1$.
\end{enumerate}

\begin{lemma}
 Let $(\varepsilon_{n},n\in\mathbb{Z})$  be an i.i.d. sequence of  $\textbf{B}-$valued random variables belongs to
$\textbf{\textrm{L}}_\textbf{B}^{p}$, $p\geq1$. Assume the conditions $(C1(p))-(C3(p))$ hold.
Then the equation \eqref{eq:autoregressive equation} has a unique strictly stationary solution given by
\begin{equation}
X_{n}=\mu+\sum \limits_{j =0 }^\infty{A_{n,j}\varepsilon_{n-j}}, \ \ n\in\mathbb{Z},\label{eq:solution of BRCA}
\end{equation}
where $A_{n,0}=I$ and $A_{n,j}=\rho_{n}\ast...\ast\rho_{n-j+1}$ for $j\geq1 $, the series on the right-hand side of \eqref{eq:solution of BRCA} converges absolutely almost surely and in $\textbf{\textrm{L}}_\textbf{B}^{p}$, $p\geq1$.
\end{lemma}
\begin{remark} 1) For finite-dimensional case $ \mathbb{R}^{d}$, $d\geq1$, conditions ensuring  the strictly stationary solution of equation \eqref{eq:autoregressive equation} are well-known. For instance, for real valued case, we refer to \cite{paper 4}, \cite{paper 11}, \cite{book 4}, \cite{paper 6} and for multivariate case, see \cite{paper 11}, \cite{paper 9}. For real valued autoregressive process with random coefficients (RCA) the strictly stationary solution was characterized under optimal conditions. Indeed, Aue et al.\cite{paper 6} proved if $\textbf{E}[\ln^{+}|\rho_{0}|]<0 $  and $\textbf{E}[\ln^{+}|\varepsilon_{0}|]<0$  are finite, where $\ln^{+}x:=\max\{0,x\} $ is a positive part of natural logarithm, then the following condition:
\begin{equation}
-\infty\leq\textbf{E}[\ln|\rho_{0}|]<0\label{eq:condition of operator for real case}
\end{equation}
is necessary and sufficient for the existence and uniqueness of the strictly stationary solution of RCA model. It is worth noting that without involving any moment conditions RCA has a strictly stationary solution. We refer to \cite{paper 9}, \cite{book 4} for multivariate case.

2) Recently, Boukhair and Mourid \cite{paper 10} considered Hilbert space valued random coefficients autoregressive process of order 1 (HRCA(1)), and studied the existence and uniqueness of  strictly stationary solution of \eqref{eq:autoregressive equation} under conditions: $(C1), (C2)$ and
\begin{equation}
\textbf{E}[\ln\|\rho_{0}\|_{\mathcal{L}}]<0. \label{eq:condition of operator for Banach space}
\end{equation}
The condition \eqref{eq:condition of operator for Banach space} is an extension of \eqref{eq:condition of operator for real case} to Banach space setting and these conditions are also best assumptions to ensure the existence and uniqueness of  strictly stationary solution of BRCA(1) process. Clearly, condition $(C3(p))$ is stronger then \eqref{eq:condition of operator for Banach space}, however, it is a flexible for our framework. We only consider $(C3(p))$ in order to deal with some moment assumptions on process defined by \eqref{eq:autoregressive equation}. On the other hand, Boukhair and Mourid \cite{paper 10} replaced the condition $(C3(p))$ by the following condition:
\begin{equation}
 \Delta:=\sup\{\|\rho_{n}\|_{\mathcal{L}}, n\in\mathbb{Z}\}<1,\ \ a.s.\label{eq:condition of Mourid}
\end{equation}
Obviously, $(C3(p))$ is weaker than \eqref{eq:condition of Mourid}.
\end{remark}
\begin{lemma} Let $(X_{n},n\in\mathbb{Z})$  be BRCA(1) process with conditions $(C1(2))-(C3(2))$. Then
\begin{equation}
C_{X_{0}}=\textbf{E}(\rho_{1}C_{X_{0}}\rho_{1}^{*})+C_{\varepsilon_{0}},\label{eq:covariance operator of BRCA}
\end{equation}
\end{lemma}
The following result states a weak law of large numbers for BRCA(1) process in nuclear norm $ \| \cdot\|_{\mathcal{N}} $.
\begin{theorem} Let $(X_{n},n\in\mathbb{Z})$ be BRCA(1) process belongs to $\textbf{\textrm{L}}_\textbf{B}^{2}$ space. If conditions $(C1(2))-(C3(2))$ hold, then
\begin{equation}
\parallel nC_{\overline{X}_{n}}-\sum \limits_{h=-\infty}^\infty{C_{X_{0},X_{h}}}\parallel_{\mathcal{N}} \rightarrow 0, \ \ n\rightarrow\infty.\label{eq:Strong law}
\end{equation}
\end{theorem}
We will give one immediate corollary of Theorem 3.4.

\begin{corollary} Let $\textbf{B}=\textbf{H}$ be a Hilbert space with inner product $\langle\cdot,\cdot\rangle_{\textbf{H}}$ which generates norm $\|\cdot\|_{\textbf{H}}$. Consider HRCA(1) process with conditions $(C1(2))-(C3(2))$. Then, we may give the exact rate of convergence to \eqref{eq:Strong law}.
\end{corollary}
Indeed, in this case,
\begin{center}
$n\textbf{E}\| \overline{X}_{n}-\mu\ \|_{\mathbf{H}}^{2}\rightarrow \sum \limits_{h=-\infty}^\infty{\textbf{E}\langle X_{0},X_{h}\rangle_{\mathbf{H}}}, \ \ n\rightarrow\infty.$
\end{center}

\begin{remark} Note that due to the inequality $\left\| C_{X} \right\|_{{\rm {\mathcal N}}} \le \textbf{E} \left\| X\right\| _{\textbf{B} }^{2}$ in Banach space setting, one can not directly obtain asymptotic properties of $\textbf{E}\left\| \frac{S_{n} }{n} \right\|^{2}$ from \eqref{eq:Strong law}.
\end{remark}
The next result concerns consistency of the sample mean ${\overline{X}_{n}}$.

\begin{theorem} Let $(X_{n},n\in\mathbb{Z})$ be BRCA(1) process defined in $\textbf{\textrm{L}}_\textbf{B}^{1}$ with $\textbf{E} X_{0}=0$. If conditions $(C1(1))-(C3(1))$ hold, then $\overline{X}_{n}$ converges to $\mu$ almost surely.
\end{theorem}
The following result deals with complete convergence for BRCA(1) process.

\begin{theorem} Let $(X_{n},n\in\mathbb{Z})$ be BRCA(1) process satisfying conditions $(C1(p))-(C3(p))$. Assume that $(\textbf{B}, \|\cdot\|_{\textbf{B}}) $ is a $r-$smooth Banach space for some $r>p$. If ${\rm E}\left\| \varepsilon_{0} \right\|_{{\mathbf{B}}}^{p}<\infty $ for some $p\in \left(1,2\right)$, then, for any $1\le {1\mathord{\left/ {\vphantom {1 \alpha }} \right. \kern-\nulldelimiterspace} \alpha } \le p$  and for all $\varepsilon>0$,
\begin{equation}
\sum_{n=1}^{\infty}n^{\alpha p-2} {\mathbf{P}}\left({\mathop{\max }\limits_{1\le k\le n}} \left\| S_{k} \right\|_{\textbf{B}} \ge \varepsilon n^{\alpha}\right) <\infty.\label{eq:rate of convergence}
\end{equation}
\end{theorem}
\begin{remark} Note that relation \eqref{eq:rate of convergence} describes speed of convergence in strong law of large numbers for $r-$smooth Banach space. Indeed, the sequence ${\max\limits_{1\le k\le n}} \left\| S_{k} \right\| _{\textbf{B}}$ is a monotonic and \eqref{eq:rate of convergence} holds if and only if
\begin{center}
$ \sum\limits_{n=1}^{\infty}{{\mathbf{P}}\left({\max\limits_{1\le k\le 2^{n}}} \left\| S_{k} \right\|_{\textbf{B}}\geq\varepsilon 2^{{n \mathord{\left/ {\vphantom {1 2}} \right.
\kern-\nulldelimiterspace} p}}\right)<\infty},$
\end{center}
with $ \alpha p=1$. Therefore, we may claim  that $n^{{1 \mathord{\left/ {\vphantom {1 2}} \right.
\kern-\nulldelimiterspace} p}}S_{n}\rightarrow 0 $ a.s.
\end{remark}
Now we deal with some exponential inequalities for equation \eqref{eq:autoregressive equation}. We will replace the condition $(C3(p))$ by \eqref{eq:condition of Mourid}. We say that BRCA(1) process $X_{0}$  belongs to the class $\mathcal{E}$ if for some $\gamma >0$, ${\textbf{E}}\exp\left(\gamma \left\| X_{0} \right\|_{{ \mathbf{B}}}\right)<\infty$. The following Lemma asserts on necessary and sufficient condition ensuring $X_{0} \in {\rm {\mathcal E}}$.
\begin{lemma}
Assume for BRCA(1) process conditions $(C1(p))$, $(C2)$ and \eqref{eq:condition of Mourid} hold. Then relations $X_{0} \in {\rm {\mathcal E}}$ and $\varepsilon_{0} \in {\rm {\mathcal E}}$  are equivalent.
\end{lemma}
We now deal with CLT for BRCA(1) process.

\begin{theorem} Let $(\textbf{B}, \|\cdot\|_{\textbf{B}})$  be 2-smooth Banach space. Assume the conditions $(C1(2))-(C3(2))$ hold for BRCA(1) process. Then
\begin{center}
$\sqrt{n}\left(\frac{S_{n} }{n}-\mu \right){\mathop{\to}\limits^{{\rm {\mathcal D}}}} N\sim {\rm {\mathcal N}}\left(0,\Gamma \right), \ \ n\to \infty,$
\end{center}
where $\Gamma=\left(I-{\textbf{E}}\left(\rho_{0} \right)\right)^{-1}\left(X_{1}-{\textbf{E}}\left(\rho_{0} \right)X_{0}\right)$.
\end{theorem}

\section{Proofs}
\textbf{Proof of Lemma 3.1.} By Minkowski inequality and from conditions $(C1(p))-(C3(p))$, we have as $k,l\rightarrow\infty $
\begin{center}
$ \textbf{E} \left\|\sum\limits_{j=k}^{l}{A_{n,j} \varepsilon_{n-j}}  \right\|_{\textbf{B} }^{p} \leq \left[\sum \limits_{j=k}^{l}\left(\textbf{E} \left(\left\| A_{n,j} \right\|_{{\rm {\mathcal L}}\left({\rm B} \right)}^{p} \right)\right)^{{1\mathord{\left/ {\vphantom {1 p}} \right. \kern-\nulldelimiterspace} p} } \left({\rm E} \left(\left\| \varepsilon _{n-j} \right\|_{\textbf{B} }^{p} \right)\right)^{{1\mathord{\left/ {\vphantom {1 p}} \right. \kern-\nulldelimiterspace} p} }  \right]^{p}$
\end{center}
\begin{center}
$=\textbf{E} \left(\left\| \varepsilon _{0} \right\| _{\textbf{B} }^{p} \right)\left[\sum\limits _{j=k}^{l}\left(\textbf{E} \left\| A_{n,j} \right\|_{{\rm {\mathcal L}}\left(\textbf{B} \right)}^{p} \right)^{{1\mathord{\left/ {\vphantom {1 p}} \right. \kern-\nulldelimiterspace} p} }\right]^{p}  \leq  \textbf{E} \left(\left\| \varepsilon _{0} \right\| _{\textbf{B} }^{p} \right)\left[\sum \limits_{j=k}^{l}\left(\textbf {E} \left\| \rho_{0} \right\|_{{\rm {\mathcal L}}}^{p} \right)^{{j\mathord{\left/ {\vphantom {j p}} \right. \kern-\nulldelimiterspace} p} }\right]^{p}  \to 0. $
\end{center}
Thus from the Cauchy criterion we may deduce that the series in \eqref{eq:solution of BRCA} converges in $\textbf{\textrm{L}}_\textbf{B}^{1}$.\\
It remains to prove the almost sure convergence of the series in \eqref{eq:solution of BRCA}.

Observe that
\begin{center}
$\textbf{E} \left(\sum\limits_{j=0}^{\infty }\left\| A_{n,j} \varepsilon_{n-j} \right\|_{\mathcal{L(\textbf{B})}}  \right)^{p} ={\mathop{\lim }\limits_{n\to \infty }} \textbf{E} \left(\sum\limits_{j=0}^{n}\left\| A_{n,j} \varepsilon_{n-j} \right\|_{\mathcal{L(\textbf{B})}}  \right)^{p}$
\end{center}
\begin{center}
$\leq {\mathop{\lim}\limits_{n\to \infty }} \left[\sum\limits_{j=1}^{n}\left(\textbf{E} \left(\left\| A_{n,j} \right\|_{\mathcal{L}} \left\| \varepsilon_{n-j} \right\| \right)^{p} \right)^{{1\mathord{\left/ {\vphantom {1 p}} \right. \kern-\nulldelimiterspace} p} } \right]^{p} \leq \left(\textbf{E} \left\| \varepsilon _{0} \right\|_{\textbf{B} }^{p} \right)\left[\sum\limits_{j=0}^{\infty }\left(\textbf{E} \left\| \rho _{0} \right\|_{{\rm {\mathcal L}}\left(\textbf{B} \right)}^{p} \right)^{{j\mathord{\left/ {\vphantom {j p}} \right. \kern-\nulldelimiterspace} p}} \right]^{p} <\infty$
\end{center}
and therefore the series $\sum\limits_{j=0}^{\infty }A_{n,j}\varepsilon_{n-j}$ converges absolutely and almost surely.

Let us consider the stationary process
\begin{center}
$Y_{n} =\sum\limits_{j=0}^{\infty}A_{n,j} \varepsilon_{n-j}, \ \ n\in {\mathbb{ Z}}.$
\end{center}
Now we check that $\left(Y_{n}, n\in {\rm \mathbb{Z}}\right)$  is a solution of \eqref{eq:autoregressive equation}. Then iterating, we get
\begin{center}
$Y_{n}-\rho_{n} Y_{n-1}=\sum\limits_{j=0}^{\infty }A_{n,j} \varepsilon_{n-j}-\sum\limits_{j=0}^{\infty}\rho_{n} A_{n-1,j} \varepsilon _{n-1-j}$
\end{center}
\begin{center}
$ =\sum\limits_{j=0 }^\infty{A_{n,j} \varepsilon_{n-j}}-\sum\limits_{j=0}^\infty{A_{n,j+1} \varepsilon_{n-1-j}}$
\end{center}
\begin{center}
$ =\sum\limits_{j=0 }^\infty{A_{n,j} \varepsilon_{n-j}}-\sum\limits_{i=1}^\infty{A_{n,i} \varepsilon_{n-i}}=\varepsilon _{n}$
\end{center}
Hence, $\left(Y_{n}, n\in {\mathbb{Z}}\right)$ is a stationary solution of BRCA(1).

Conversely, let us suppose that $(X_{n},n\in \mathbb{Z})$ is a another stationary solution of equation \eqref{eq:autoregressive equation}. Then
from the decomposition of $\left(X_{n} ,n\in \mathbb{Z}\right)$, it follows
\begin{equation}
X_{n}=\sum\limits_{j=0}^{n-1}A_{n,j} \varepsilon_{n-j}+A_{n,n}  X_{0},\ \ n\geq1. \label{eq:finite decomposition}
\end{equation}
Now taking into account the stationarity of $\left(X_{n}, n\in \mathbb{Z}\right)$ and by condition $(C3(p))$, one has
\begin{center}
$\left\| X_{n}-\sum\limits_{j=0}^{n-1}A_{n,j} \varepsilon_{n-j}  \right\|_{{\rm L}_{\textbf{B} }^{p} }^{p} \leq \left\| \textbf{E} \left(A_{n,n} X_{0} \right)\right\|_{\textbf{B} }^{p}=\textbf{E} \left(\left\| X_{0} \right\|_{\textbf{B} }^{p} \right)\left(\left\| \textbf{E} \rho_{0} \right\|_{{\rm {\mathcal L}}\left(\textbf{B} \right)}^{p} \right)^{n} \to 0,\ \  n\to \infty.$
\end{center}
which proves uniqueness of stationary solution of BRCA(1). This ends the proof of Lemma.

\textbf{Proof of Lemma 3.3.} From equation \eqref{eq:autoregressive equation} we obtain
\begin{center}
$C_{X_{0}}({x^{*}})=C_{X_{1}}({x^{*}})=\textbf{E}[x^{*}(X_{1})X_{1}]=\textbf{E}[x^{*}(\rho_{1}X_{0})+\varepsilon_{1})\rho_{1}X_{0}+\varepsilon_{1}]$
\end{center}
\begin{center}
$=\textbf{E}(x^{*}(\rho_{1}X_{0})\rho_{1}X_{0})+\textbf{E}(x^{*}(\rho_{1}X_{0})\varepsilon_{1})+\textbf{E}(x^{*}(\varepsilon_{1})\rho_{1}X_{0})+ \textbf{E}(x^{*}(\varepsilon_{1})\varepsilon_{1}).$
\end{center}
By independence $\rho_{1}X_{0}$ of $\varepsilon_{1}$, we have
\begin{center}
$\textbf{E}(x^{*}(\rho_{1}X_{0})\varepsilon_{1})=\textbf{E}(x^{*}(\varepsilon_{1})\rho_{1}X_{0})=0$
\end{center}
which implies \eqref{eq:covariance operator of BRCA}.

\textbf{Proof of Theorem 3.4.} By independence $\left(A_{h,j} \varepsilon _{h-j},0\le j\le h-1,h\ge 1\right)$ of $X_{0}$, we deduce
\begin{equation}
C_{X_{0}, X_{h}}=C_{X_{0},A_{h,h} \left(X_{0} \right)}.\label{eq:cross-covariance}
\end{equation}
From \eqref{eq:cross-covariance} it follows that
\begin{equation}
\left\| C_{X_{0},X_{h} } \right\|_{{\rm {\mathcal N}}} \leq \left(\textbf {E} \left\| \rho_{0} \right\|_{{\rm {\mathcal L}}\left(\textbf{B} \right)} \right)^{h} \textbf {E} \left\| X_{0} \right\|_{\textbf {B} }^{2},\ \  h\geq 1.\label{inequality for cross-covariance}
\end{equation}
Therefore, we may deduce that the series $\sum\limits_{h=-\infty}^{\infty}C_{X_{0}, X_{h}}$ converges in the space of nuclear operators over $\textbf{B}$.

 On the other hand, one has
\begin{center}
$y^{*} \left(C_{\overline{X}_{n}} \left(x^{*} \right)\right)=\textbf{E} \left(x^{*} \left(\overline{X}_{n}\right)y^{*} \left(\overline{X}_{n} \right)\right)=\frac{1}{n^{2} } \sum\limits_{1\leq i,j\leq n}\textbf{E}  \left(x^{*} \left(X_{i} \right)y^{*} \left(X_{j} \right)\right), x^{*},y^{*} \in \textbf{B}^{*}.$
\end{center}
By stationarity of $\left(X_{n}, n\in {\mathbb{Z}}\right)$, we get
\begin{center}
$y^{*} \left(nC_{\overline{X}_{n}} \left(x^{*} \right)\right)=\sum\limits_{\left|h\right|\le n-1}\left(1-\frac{\left|h\right|}{n} \right) \textbf{E} x^{*} \left(X_{0} \right)y^{*} \left(X_{h} \right),$
\end{center}
and it yields
\begin{center}
$nC_{\overline{X}_{n}}=\sum\limits_{\left|h\right|\leq n-1}\left(1-\frac{\left|h\right|}{n} \right) C_{X_{0}, X_{h} }.$
\end{center}
Hence it remains to note that \eqref{eq:Strong law} is a consequence of \eqref{inequality for cross-covariance}.

\textbf{Proof of Theorem 3.7.} The fact that $\left(\varepsilon _{n}, n\in {\mathbb{Z}}\right)$ is a ${\textbf{B}} -$valued white noise and $\left(\rho_{n}, n\in {\mathbb{Z}}\right)$ is an i.i.d. sequence ensure the ergodicity of the sequence$\left(X_{n},n\in {\mathbb{Z}}\right)$. Hence by Mourier's ergodic theorem \cite{paper 31} for stationary Banach-valued random variables it follows $\overline{X}_{n}$ converges to $\mu$ almost surely.

\textbf{Proof of Theorem 3.8.} Clearly,
\begin{center}
${\mathop{\max}\limits_{1\leq k\leq n}} \left\| S_{k} \right\|_{\textbf{B} } \leq \left\| N_{0} \right\|_{\textbf{B} }+{\mathop{\max }\limits_{2\le k\leq n}} \left\| N_{k} \right\| _{\textbf{B} }+{\mathop{\max }\limits_{1\leq k\leq n}} \left\| \sum\limits_{i=1}^{k}M_{i}  \right\| _{\textbf{B}}.$
\end{center}
From Woyczy\'{n}ski \cite{paper 37}, we know that
\begin{center}
$\sum\limits_{n=1}^{\infty }n^{\alpha p-2} \textbf{P} \left({\mathop{\max}\limits_{1\leq k\leq n}} \left\| \sum\limits _{i=1}^{k}M_{i}  \right\|_{\textbf{B}} \geq \varepsilon n^{\alpha } \right) <\infty.$
\end{center}
Hence it remains to verify that
\begin{center}
$\sum\limits_{n=1}^{\infty }n^{\alpha p-2} \textbf{P} \left({\mathop{\max }\limits_{2\le k\le n}} \left\| N_{k} \right\| _{\textbf{B}} \geq \varepsilon n^{\alpha } \right) <\infty.$
\end{center}
Since $\left(X_{n}, n\in {\mathbb{Z}}\right)$ is a strictly stationary, we obtain
\begin{center}
$\sum\limits_{n=1}^{\infty}n^{\alpha p-2} \textbf{P} \left({\mathop{\max }\limits_{2\le k\le n}} \left\| N_{k} \right\| _{\textbf{B}} \geq \varepsilon n^{\alpha } \right)$
\end{center}
\begin{center}
$\leq\sum\limits_{n=1}^{\infty }n^{\alpha p-1} \textbf{P} \left(\left\| X_{0} \right\|_{\textbf{B} } \geq \varepsilon n^{\alpha } \left\| \left(I-\bar{\rho }\right)^{-1} \right\| _{{\rm {\mathcal L}}\left(\textbf{B} \right)}^{-1} \left\| \bar{\rho }\right\| _{{\rm {\mathcal L}}\left(\textbf{B} \right)}^{-1} \right).$
\end{center}
Applying Fubini's theorem and noting that $X_{0}$ belongs to $\textbf{\textrm{L}}_\textbf{B}^{p}$, $p\geq1$, we deduce that the last series is finite.

\textbf{Proof of Lemma 3.10.} Let us suppose that $X_{0}\in \mathcal{E}$. Then from autoregressive equation \eqref{eq:autoregressive equation} and \eqref{eq:condition of Mourid} we have
\begin{center}
$\|\varepsilon_{n}\|_{\textbf{B}}\leq\|\rho_{n}\|_{\mathcal{L}}\|X_{n-1}\|_{\textbf{B}}+\|X_{n}\|_{\textbf{B}}\leq  \Delta\|X_{n-1}\|_{\textbf{B}}+\|X_{n}\|_{\textbf{B}}. \ \ $
\end{center}
Hence for any $\alpha>0 $, a.s.
\begin{center}
$ \textbf{E}(\exp \alpha \|\varepsilon_{1}\|_{\textbf{B}})\leq \textbf{E}(\exp (\alpha \Delta\|X_{0}\|_{\textbf{B}}+\alpha\|X_{1}\|_{\textbf{B}})$
\end{center}
\begin{center}
$\leq \textbf{E}\exp (2\alpha \|X_{0}\|_{\textbf{B}}).$
\end{center}
Now it remains to choose $\alpha$, and if we choose $\alpha=\frac{\gamma}{2}$ and taking into account the assumption $X_{0}\in \mathcal{E}$ we obtain
\begin{center}
$ \textbf{E}(\exp \alpha \|\varepsilon_{0}\|_{\textbf{B}})\leq \textbf{E}(\exp \gamma \|X_{0}\|_{\textbf{B}})<\infty.$
\end{center}
Therefore, we get $\varepsilon_{0}\in \mathcal{E}$.

Conversely, assume that $X_{0}\in \mathcal{E}$ for some $\gamma>0$. From \eqref{eq:condition of Mourid} it follows that there exists $c>0$ such that $\Delta\leq c<1$ a.s. Therefore, from \eqref{eq:solution of BRCA} it entails that
\begin{center}
$ \|X_{n}\|_{\textbf{B}}\leq \sum\limits_{j=0}^{\infty }{\|A_{n,j}\|_{\mathcal{L}}\|\varepsilon_{n-j}\|_{\textbf{B}}}\leq\sum\limits_{j=0}^{\infty }{\Delta^{j}\|\varepsilon_{n-j}\|_{\textbf{B}}}.$
\end{center}
Due to independence of $\|\varepsilon_{n}\|_{\textbf{B}}, n\in \mathbb{Z}$ and using monotone convergence of expectation, we obtain for all $\alpha>0$
\begin{center}
$ \textbf{E}(\exp \alpha \|X_{n}\|_{\textbf{B}})\leq \prod\limits_{j=0}^{\infty}{\textbf{E}(\exp \alpha c^{j}\|\varepsilon_{n-j}\|_{\textbf{B}})}.$
\end{center}
Now, in order to verify the convergence of last infinite product, it suffices to prove
\begin{center}
$\sum\limits_{j=0}^{\infty}{\log\textbf{E}(\exp \alpha c^{j}\|\varepsilon_{n-j}\|_{\textbf{B}})}<\infty. $
\end{center}
Notice that
\begin{center}
$\lim\limits_{n\rightarrow\infty}\sum\limits_{j=0}^{n}{\log\textbf{E}(\exp \alpha c^{j}\|\varepsilon_{n-j}\|_{\textbf{B}})}\leq \lim\limits_{n\rightarrow\infty}\sum\limits_{j=0}^{n}{\log(1+\sum\limits_{i=1}^{\infty}{\frac{\alpha^{i}(c^{j})^{i}}{i!}\textbf{E}\|\varepsilon_{0}\|_{\textbf{B}}^{i}}) }$
\end{center}
\begin{center}
$\leq \sum\limits_{i=1}^{\infty}{\frac{\alpha^{i}(c^{j})^{i}}{i!}\textbf{E}\|\varepsilon_{0}\|_{\textbf{B}}^{i}} \lim\limits_{n\rightarrow\infty}\sum\limits_{j=0}^{n}{(c^{j})^{i}}\leq \textbf{E}\exp(\gamma \|\varepsilon_{0}\|_{\textbf{B}}) \sum\limits_{j=0}^{\infty}{(c^{j})}<\infty.$
\end{center}
Consequently, we obtained $X_{0}\in \mathcal{E}$ which ends the proof of Lemma.

\textbf{Proof of Theorem 3.11.} We again use the decomposition \eqref{eq:sum of autoregressive process} and write
\begin{equation}
\frac{S_{n}}{\sqrt{n}}=\frac{1}{\sqrt{n} } \left(N_{0}-N_{n} \right)+\frac{1}{\sqrt{n}} \sum_{i=1}^{n}M_{i}.\label{eq:decomposition for CLT}
\end{equation}
We show the convergence of the first term of \eqref{eq:decomposition for CLT} to zero in probability. Indeed, since $X_{n} $ is a strictly stationary
\begin{equation}
\frac{1}{\sqrt{n} } \left\| \left(N_{0}-N_{n} \right)\right\| \le \frac{\left\| \left(I-\bar{\rho }\right)^{-1} \right\| \left\| X_{0} \right\|_{\textbf {B} } }{\sqrt{n} } \to 0, \ \ n\to \infty.
\end{equation}
For the second term, we may apply Woyczy\'{n}ski's \cite{paper 38} CLT for $2-$smooth Banach-valued martingale differences and obtain
\begin{equation}\label{weak convergence to normal distribution}
\frac{1}{\sqrt{n} } \sum\limits_{i=1}^{n}M_{i}  {\mathop{\to }\limits^{{\rm {\mathcal D}}}} N\sim {\rm {\mathcal N}}\left(0,\left(I-\textbf{E} \left(\rho _{0} \right)\right)^{-1} \left(X_{1} -\textbf {E} \left(\rho _{0} \right)X_{0} \right)\right).
\end{equation}
Combining relations \eqref{eq:decomposition for CLT}-\eqref{weak convergence to normal distribution} we get the desired result.

Now we provide an example of BRCA(1) process.

\begin{example} Let ${\rm L}_{{\rm B} }^{p} \left(\left[0,1\right],{\rm {\tt B}}_{\left[0,1\right]} ,\lambda \right)$, $p\ge 1$ and consider a sequence of random kernel operators $\rho_{n} :{\rm L}_{{\rm B}}^{p} \to {\rm L}_{{\rm B}}^{p}$ defined by
\begin{center}
$\rho_{n} \left(x\right)\left(t\right)=\int\limits_{0}^{1}K_{n} \left(t,s\right) x\left(s\right)ds, \ \ t\in \left[0,1\right],\ \  x\in{\textbf{B}},\ \ n\in {\mathbb{Z}},$
\end{center}
where $\left(K_{n}, n\in {\mathbb{Z}}\right)$ is a sequence of i.i.d. kernels in ${\rm L}_{{\rm B}}^{p} \left(\left[0,1\right]\right)$ satisfying
\begin{center}
${\rm E} \int\limits_{0}^{1}\int\limits_{0}^{1}K_{0}^{p} \left(t,s\right)dtds<1.$
\end{center}
Take an i.i.d. mean zero innovations $\varepsilon_{n} \in {\rm L}_{{\rm B}}^{p} \left(\left[0,1\right]\right)$, independent of $(\rho _{n},n\in \mathbb{Z})$ and set
\begin{center}
${\rm E} \left\| \varepsilon_{n} \right\|_{{\rm B}}^{p}={\rm E} \left[\int\limits_{0}^{1}\varepsilon_{n}^{p}\left(t\right)dt \right]<\infty.$
\end{center}
Then the equation \eqref{eq:autoregressive equation} can be written more explicitly as
\begin{center}
$X_{n} \left(t\right)=\int\limits_{0}^{1}K_{n} \left(t,s\right) X_{n-1} \left(s\right)ds+\varepsilon _{n} \left(t\right), \ \ t,s\in \left[0,1\right].$
\end{center}
\end{example}


\end{document}